\documentclass[12pt]{amsart}

\usepackage{geometry}
\usepackage{hyperref}

\renewcommand{\phi}{\varphi}

\newcommand{\zz}{\mathbb{Z}}
\newcommand{\cc}{\mathbb{C}}

\newcommand{\oo}{\mathcal{O}}

\newtheorem{thm}{Theorem}[section]
\newtheorem{lm}[thm]{Lemma}
\newtheorem*{conj}{Conjecture}
\newtheorem{corollary}[thm]{Corollary}

\theoremstyle{definition}
\newtheorem{defi}{Definition}

\theoremstyle{remark}
\newtheorem*{rem}{Remark}

\begin{document}

\title{ Integrality Properties of the CM-values of Certain Weak Maass Forms}

\author{Eric Larson}
\address{Department of Mathematics. Harvard University, Cambridge, MA 02138.}
\email{elarson3@gmail.com}

\author{Larry Rolen}
\address{Department of Mathematics and Computer Science, Emory University, Atlanta, GA 30322.}
\email{larry.rolen@mathcs.emory.edu}

\begin{abstract}
In a recent paper, Bruinier and Ono prove that the coefficients of certain weight $-1/2$ harmonic Maass forms are traces of singular moduli for weak Maass forms. In particular, for the partition function $p(n)$, they prove that \[p(n)=\frac{1}{24n-1}\cdot\displaystyle\sum P(\alpha_Q),\] where $P$ is a weak Maass form and $\alpha_Q$ ranges over a finite set of discriminant $-24n+1$ CM points. Moreover, they show that $6\cdot(24n-1)\cdot P(\alpha_Q)$ is always an algebraic integer, and they conjecture that $(24n-1)\cdot P(\alpha_Q)$ is always an algebraic integer. Here we prove a general theorem which implies this conjecture as a corollary. 
\end{abstract}
\maketitle

\section{Introduction and Statement of Results \label{sec:intro}}
A \emph{partition} of a positive integer $n$ is any nonincreasing sequence of  positive integers which sum to $n$. The partition function $p(n)$, which counts the number of partitions of $n$, is an important function in number theory whose study has a long history. One of the celebrated results of Hardy and Ramanujan on this function, giving rise to the ``circle'' method, quantifies the growth rate: \[p(n)\sim\frac{1}{4n\sqrt{3}}\cdot e^{\pi\sqrt{2n/3}}.\] This asymptotic and its method of proof were later refined by Rademacher, yielding an ``exact'' formula in terms of a modified Bessel function of the first kind $I_{\frac{3}{2}}(\cdot)$ and a Kloosterman sum $A_k(n)$: \[p(n)=2\pi(24n-1)^{-\frac{3}{4}}\displaystyle\sum_{k=1}^{\infty}\frac{A_k(n)}{k}\cdot I_{\frac{3}{2}}\left(\frac{\pi\sqrt{24n-1}}{6k}\right).\] One can compute values of $p(n)$ from this formula by using sufficiently accurate truncations. Bounding the resulting error is a well-known difficult problem; the best-known bounds are due to Folsom and Masri \cite{Folsom-Masri}. 

In recent work \cite{Ono-p(n)}, Bruinier and Ono prove a new formula for $p(n)$ as a finite sum of algebraic numbers. These numbers are \emph{singular moduli} for a \emph{weak Maass form} which they describe in terms of Dedekind's eta function and the quasimodular Eisenstein series $E_2$,
which are defined in terms of $q := e^{2 \pi i z}$ as
\[\eta(z):=q^{\frac{1}{24}} \prod_{n=1}^{\infty}(1-q^n) \quad \text{and} \quad E_2(z):=1-24\displaystyle\sum_{n=1}^{\infty}\displaystyle\sum_{d|n}dq^n.\]
They then define the $\Gamma_0(6)$ weight $-2$ meromorphic modular form:
\begin{equation} F_p(z):=\frac{1}{2}\cdot\frac{E_2(z)-2E_2(2z)-3E_2(3z)+6E_2(6z)}{\eta(z)^2\eta(2z)^2\eta(3z)^2\eta(6z)^2}=q^{-1}-10-29q-\cdots\end{equation}
Using the convention that $z:=x+iy$, with $x,y\in\mathbb{R}$, they define the weak Maass form:
\begin{equation}P_p(z):=-\left(\frac{1}{2\pi i}\cdot \frac{d}{dz}+\frac{1}{2\pi y}\right)F_p(z)=\left(1-\frac{1}{2\pi y}\right)q^{-1}+\frac{5}{\pi y}+\left(29+\frac{29}{2\pi y}q\right)+\cdots\end{equation}
Bruinier and Ono give a formula for $p(n)$ in terms of
discriminant $-24n+1=b^2-4ac$ positive definite integral binary quadratic forms $Q(x,y)=ax^2+bxy+cy^2$ satisfying the condition $6\mid a$. The group $\Gamma_0(6)$ acts on such forms, and we let $\mathcal{Q}_n$ be any set of representatives of those equivalence classes with $a>0$ and $b\equiv 1\ (\text{mod } 12)$.
To each such $Q$, we associate the CM point $\alpha_Q$ defined to be the root of $Q(x,1)=0$ lying in the upper half of the complex plane. Then the formula of Bruinier and Ono states: \begin{equation}p(n)=\frac{1}{24n-1}\cdot\displaystyle\sum_{Q\in\mathcal{Q}_n }P(\alpha_Q).\end{equation}
They further prove that each $6\cdot(24n-1)\cdot P_p(\alpha_Q)$ is an algebraic integer. They also show that the numbers $P(\alpha_Q)$, as $Q$ varies over $\mathcal{Q}_n$, form a multiset which is a union of Galois orbits for the discriminant $-24n+1$ ring class field. Based on numerics, they made the following conjecture:
\begin{conj}[Bruinier and Ono \cite{Ono-p(n)}]\label{conjecture} For the Maass form $P(z)$ above and for the $\alpha_Q$ in the formula for $p(n)$, we have that $(24n-1)\cdot P_p(\alpha_Q)$ is an algebraic integer.
\end{conj}
We prove that this is indeed the case. In fact, we prove that this is the true for all the CM points of discriminant $-24n+1$ of a wider class of Maass forms. Namely, we have the following:
\begin{thm}\label{mainthm}
Suppose $F$ is a weakly holomorphic, weight $-2$ modular form on a congruence subgroup
such that the Fourier expansions of
\[F \quad \text{and} \quad q\frac{dF}{dq} + F \cdot \frac{E_2 E_4 - E_6}{6E_4}\]
at all
cusps have coefficients that are algebraic integers.
Let $\alpha_Q$ be the CM point in $\mathbb{H}$ corresponding to a
quadratic form $Q(x,y)$ of discriminant $-24n+1$,
and let $P(z)$ be the weak Maass form
\[P(z) = -\left(\frac{1}{2\pi i}\cdot\frac{d}{dz}+\frac{1}{2\pi y}\right)F(z).\]
Then $(24n-1)\cdot P(\alpha_Q)$ is an algebraic integer. 
\end{thm}
\begin{rem} We recall that a meromorphic modular form is said to be \emph{weakly holomorphic} if its poles are supported on the cusps.\end{rem}
The form $F_p(z)$ studied by Bruinier and Ono satisfies these conditions.
One can see this because $F_p(z)$ has level $6$, so the group of Atkin-Lehner
involutions acts transitively on the cusps. Since $F_p(z)$ is an eigenform
for all of the Atkin-Lehner involutions and has an integral Fourier
expansion at infinity, it follows that the Fourier
expansions of $F_p$ at all cusps is integral.
Moreover, since the Atkin-Lehner involutions commute with the Maass
raising operator
\[R_{-2} = -4\pi q \frac{d}{dq} - \frac{2}{y},\]
the Fourier expansion of 
\[q\frac{dF}{dq} + F \cdot \frac{E_2 E_4 - E_6}{6E_4} = F \cdot \frac{\left(E_2 - \frac{3}{\pi \operatorname{Im} z}\right) E_4 - E_6}{6E_4} - \frac{1}{4\pi} R_{-2} F\]
at all cusps is integral as well. Therefore, Theorem \ref{mainthm} implies the following:

\begin{corollary}\label{corollary}
The conjecture of Bruinier and Ono is true.
\end{corollary}
\begin{rem}
Corollary \ref{corollary} is sharp for small (and possibly all) $n$.
For example, we have
\[\prod_{m=1}^3(x-P_p(\alpha_{Q_m})) = x^3-23x^2+\frac{3592}{23}x-419,\]
where $Q_m$ ranges over any choice of representatives
of $\mathcal{Q}_1$.
\end{rem}
Returning to the general case of $P(z)$ as in Theorem \ref{mainthm},
the work of Bruinier and Ono (Theorem 4.5 of \cite{Ono-p(n)}) implies
$6 \cdot (24n - 1) \cdot P(\alpha_Q)$ is
an algebraic integer. Although this theorem is stated for squarefree level
and when $F$ is an eigenfunction of the Atkin-Lehner involutions,
an inspection of the proof shows that the assumptions in the statement
of Theorem~\ref{mainthm} are also sufficient.
Thus
it suffices to show that $P(\alpha_Q)$ is integral at primes
lying over $6$. We will henceforth refer to this property as
\emph{$6$-integrality}.
For this purpose,
it is convenient to decompose $P$ as
\begin{equation}\label{decomp} P=A+B\cdot C,\end{equation} where
\begin{align}
A &= -q\frac{dF}{dq}-\frac{1}{6}FE_2+\frac{FE_6(7j-6912)}{6E_4(j-1728)}, \\
B &= \frac{FE_6j}{E_4},\\
C &= \frac{E_4}{6 E_6 j}\left(E_2-\frac{3}{\pi\operatorname{Im}z}\right) - \frac{7j-6912}{6j(j-1728)}.
\end{align}

To establish the $6$-integrality of $P(\alpha_Q)$, it suffices
to establish the $6$-integrality of each of $A(\alpha_Q)$, $B(\alpha_Q)$
and $C(\alpha_Q)$. In Section~\ref{sec:ab}, we will use methods
similar to those of \cite{Ono-p(n)} to show that
$A(\alpha_Q)$ and $B(\alpha_Q)$ are $6$-integral.
Then in Section \ref{sec:c}, we 
show that $C(\alpha_Q)$ is $6$-integral
using a description of $C$ in terms of classical modular polynomials due to Masser.

\begin{rem} For the remainder of this paper, we fix $D \equiv 1$ mod $24$ with $D<0$, and we let $\alpha_Q$ denote any CM point of discriminant $D$.\end{rem}

\section{Proof of $6$-Integrality of $A$ and $B$ \label{sec:ab}} 
In this section, we prove the $6$-integrality of $A$ and $B$ 
at the CM-points $\alpha_Q$. We begin by showing that $j(\alpha_Q)$ is a unit at $2$ and $3$. 
\begin{lm} \label{jj}
Let $p \in \{2, 3\}$ and
$E$ be an elliptic curve defined over a number field $K$
having complex
multiplication by an order in a quadratic field $F$.
If $E$ has good ordinary reduction at all primes lying
over $p$, then $j(E)$ is coprime to $p$.
\end{lm}
\begin{proof}
Assume to the contrary that $j(E)$ was not coprime
to $p$; write $\mathfrak{p}$ for a prime ideal lying
over $p$ containing $j(E)$, and write $k$ for the residue
field $\oo_K / \mathfrak{p}$.

When $p = 2$, the elliptic curve $E^2 = \cc / \zz[\omega]$
(where $\omega$ is a primitive cube root of unity)
has good supersingular reduction at $\mathfrak{p}$.
But $j(\omega) = 0$, so $E^2_{/k} \simeq E_{/k}$,
so $E_{/k}$ is supersingular, which is a contradiction.

Similarly, when $p = 3$, the elliptic curve $E^3 = \cc / \zz[i]$
has good supersingular reduction at $\mathfrak{p}$.
But $j(\omega) = 1728$, so $E^3_{/k} \simeq E_{/k}$,
so $E_{/k}$ is supersingular, which is a contradiction.
\end{proof}

By this lemma, it suffices to show that both $B$ and
$A':=A\cdot j\cdot(j-1728)$ assume integral values at
all CM-points.

\begin{lm} \label{holo}
The modular functions $A'$ and $B$ are weakly
holomorphic and have integral Fourier expansions at
all cusps.
\end{lm}

\begin{proof}
By definition, we have
\[B = F \cdot E_6 \cdot \frac{j}{E_4},\]
and by direct examination, all three of the above terms
are weakly holomorphic and have integral Fourier expansions
at all cusps.
Similarly, by definition, we have
\[A' = F \cdot E_6 (j - 864) \cdot \frac{j}{E_4} - (j - 1728) \cdot \left[j \cdot \left(q \frac{dF}{dq} + F \cdot \frac{E_2 E_4 - E_6}{6E_4}\right)\right]\]
and all of the above terms are weakly holomorphic
and have integral Fourier expansions
at all cusps.
\end{proof}

\begin{lm} \label{fourier}
A weakly holomorphic modular function $g$ for a congruence
subgroup $\Gamma_g$ that has integral
Fourier expansions at all cusps is integral at any CM-point.
\end{lm}
\begin{proof} (The following argument is due to Bruinier and Ono;
see Lemma~4.3 of \cite{Ono-p(n)}.)
We consider the polynomial
\[\Psi_g(X, z) = \prod_{\gamma \in \Gamma_g \backslash \Gamma(1)} (X - g(\gamma z)).\]
This is a monic polynomial in $X$ of degree $[\Gamma(1) : \Gamma_g]$
whose coefficients are weakly holomorphic modular functions
in $z$ for the group $\Gamma(1)$, so $\Psi_g(X, z) \in \cc[j(z), X]$.

Our assumption that $g$ has integral Fourier expansion at all cusps
implies that for any $\gamma \in \Gamma(1)$, the modular function $g \mid \gamma$
has a Fourier expansion at infinity whose coefficients are algebraic integers.
Thus, the coefficients of $\Psi_g(X, z)$ are polynomials in $j(z)$
whose coefficients are algebraic integers.

Since $j$ is integral at any CM-point $\alpha$, the value
$g(\alpha)$ satisfies a monic polynomial whose coefficients are algebraic
integers, and is therefore an algebraic integer.
\end{proof}

\begin{rem}
In Appendix~\ref{explicit}, we give
values of the polynomials $\Psi_{A'}$ and $\Psi_B$
for the form $F_p$ considered by Bruinier and Ono,
thus providing a direct proof of the integrality of $A'(\alpha_Q)$
and $B(\alpha_Q)$ in this case.
\end{rem}

\begin{lm} If $\alpha_Q$ is an CM-point with discriminant
$D \equiv 1$ mod $24$, then $A(\alpha_Q)$ and $B(\alpha_Q)$ are $6$-integral.
\end{lm}
\begin{proof}
This follows from combining Lemmas \ref{jj}, \ref{holo}, and \ref{fourier}.
\end{proof}

\section{Proof of $6$-Integrality for $C(\alpha_Q)$ \label{sec:c}} 

In this section, we finish the proof of Theorem \ref{mainthm} by showing that $C(\alpha_Q)$ is $6$-integral
at the required CM points $\alpha_Q$.
To do this, we study the classical modular polynomials
$\Phi_{-D}$, in a fashion similar to Appendix~1 of \cite{Masser}.
We begin by reviewing the definition
of $\Phi_{-D}$.

\begin{defi} We say that two matrices $B_1$ and $B_2$
are \emph{equivalent} if $B_1 = X \cdot B_2$ for some
$X \in \operatorname{SL}_2(\zz)$.
\end{defi}

It is well-known that there are only finitely many
equivalence classes of primitive integer
matrices of determinant
$-D$. Write $M_1, M_2, \ldots, M_n$ for these
equivalence classes and suppose $M_1$
is such that $\alpha_Q = M_1 \alpha_Q$.

\begin{defi} We write $\Phi_{-D}(X, Y)$
for the \emph{classical modular polynomial}, i.e.\ the
polynomial such that
\[\Phi_{-D}(j(z), Y) = \prod_{i = 1}^n (Y - j(M_i z)).\]
\end{defi}

By \cite{Serre-CM}, Theorem~1 of Section~3.4,
the polynomial $\Phi_{-D}(X, Y)$ is symmetric
in $X$ and $Y$ and has coefficients that are rational
integers. In particular, we can expand $\Phi_{-D}(X, Y)$
in a power series about $X = Y = j(\alpha_Q)$ as
\[\Phi(X, Y) = \sum_{\mu, \nu} \beta_{\mu,\nu} (X - j(\alpha_Q))^{\mu} (Y - j(\alpha_Q))^{\nu},\]
where $\beta_{\mu, \nu} = \beta_{\nu, \mu}$.
We write $\beta = \beta_{0,1} = \beta_{1,0}$.

We define $Q$ to be \emph{special} is there is more than one equivalence class of matrices $M$ such that $M\alpha_Q=\alpha_Q$. This can only happen if $D=3d^2$ for some integer $d$ (see \cite{Masser}, Appendix~1), so in particular forms of discriminant $-24n+1$ are not special. 
\begin{lm}[Masser] If $Q$ is not special, we have $\beta \neq 0$ and
\[C(\alpha_Q) = \frac{\beta_{0,2} - \beta_{1,1} + \beta_{2,0}}{\beta}.\]
\end{lm}
\begin{proof}
See \cite{Masser}, Appendix~1 (in particular, equations
(100) and (106), and the definition of $\gamma$ on page 118).
\end{proof}

By definition, the $\beta_{\mu, \nu}$
are algebraic integers. Thus, to prove that $C(\alpha_Q)$ is
integral at primes lying over $6$,
it suffices to show that $\beta$ is a unit at primes lying over $6$.
From the definition of $\beta$, we have
\[\beta = \prod_{i = 2}^n (j(\alpha_Q) - j(M_i \alpha_Q)).\]
Thus, it suffices to show that for any prime $\mathfrak{p}$
lying over $6$, we have $j(\alpha_Q) \not\equiv j(M_i \alpha_Q)$ mod $\mathfrak{p}$.
To show this, it is enough to establish the following lemma:

\begin{lm} Suppose $\mathfrak{p}$ is a prime ideal of a number field $K$. Suppose $E$ and $E^{\prime}$ are two elliptic curves over $K$ with complex multiplication by the same order $R$ in a quadratic field $F$. Suppose the index $[\mathcal{O}_F:R]$ is coprime to the residue characteristic of $\mathfrak{p}$. If both curves have good ordinary reduction at $\mathfrak{p}$ and the reduced curves are isomorphic, then $E$ and $E^{\prime}$ are also isomorphic. 
\end{lm}
\begin{proof}
Write $k$ for the residue field $\oo_K / \mathfrak{p}$
and $p = \operatorname{char}(k)$.
As the index of $R$ in $\mathcal{O}_F$ is coprime to $p$,
there is an isogeny $f\colon E\to E'$ whose degree is
coprime to $p$. Since $E$ has ordinary reduction at $\mathfrak{p}$, its endomorphisms over $k$ are a rank-2 submodule $S$ of $\mathcal{O}_F$ which contains $R$. As the index of $R$ in $\mathcal{O}_F$ is coprime to $p$, the index $d$ of $R$ in $S$ is also coprime to $p$. Choose an isomorphism between the reductions $E_{/k}$ and $E'_{/k}$. Composing this with the isogeny $f$ gives an endomorphism of $E_{/k}$, and multiplying this endomorphism by $d$ gives an endomorphism which lifts to an endomorphism $g$ of $E$ whose degree is coprime to $p$. Now the specializations of the kernels of $f\circ d$ and $g$ coincide by construction, and both kernels are subgroups whose order is coprime to $p$.
Thus, $\operatorname{ker} f \circ d = \operatorname{ker} g$,
and therefore $E\cong E'$.  
\end{proof}

This completes the proof of the $6$-integrality of $C(\alpha_Q)$,
as the assumption $D\equiv 1\mod 24$ shows that the conditions of
the above lemma are satisfied. 
By the discussion in Section~\ref{sec:intro}, this
establishes Theorem~\ref{mainthm}.

\appendix\label{explicit}

\section{The Polynomials $\Psi_{A'}$ and $\Psi_B$ for $F = F_p$ \label{explicit}}

Here, we give the explicit values of the polynomials
$\Psi_{A'}$ and $\Psi_B$ when $F = F_p$ is the form considered
by Bruinier and Ono in \cite{Ono-p(n)}. Namely, we have
\[\Psi_{A'} = X^{12} + \sum_{i = 0}^{11} a_i X^i \quad \text{and} \quad \Psi_B = X^{12} + \sum_{i = 0}^{11} b_i X^i,\]
where the $a_i$ and $b_i$ are the
polynomials in $j$ with integer coefficients given below.

\begin{align*}
a_{11} &= - 2 \cdot (j - 2^{6} \cdot 3^{3}) \cdot (j - 2^{5} \cdot 3^{3}) \cdot j \\
\displaybreak[0]
a_{10} &= - (j - 2^{6} \cdot 3^{3}) \cdot j^{2} \cdot (7 \cdot 67 \cdot j^{2} - 2^{6} \cdot 3^{2} \cdot 2053 \cdot j + 2^{11} \cdot 3^{5} \cdot 31 \cdot 53) \\
\displaybreak[0]
a_{9} &= 2 \cdot (j - 2^{6} \cdot 3^{3})^{2} \cdot j^{2} \cdot (3^{2} \cdot j^{4} - 2^{3} \cdot 6379 \cdot j^{3} + 2^{6} \cdot 3^{2} \cdot 162713 \cdot j^{2} \\
\displaybreak[0]
&\quad - 2^{12} \cdot 3^{5} \cdot 72797 \cdot j + 2^{25} \cdot 3^{12}) \\
a_{8} &= 2 \cdot (j - 2^{6} \cdot 3^{3})^{2} \cdot j^{3} \cdot (2 \cdot 7 \cdot 13^{2} \cdot j^{5} - 3^{2} \cdot 409 \cdot 3373 \cdot j^{4} \\
&\quad + 2^{7} \cdot 3^{4} \cdot 1237 \cdot 1973 \cdot j^{3} - 2^{14} \cdot 3^{7} \cdot 5 \cdot 311 \cdot 443 \cdot j^{2} \\
\displaybreak[0]
&\quad + 2^{21} \cdot 3^{10} \cdot 31 \cdot 2897 \cdot j - 2^{31} \cdot 3^{14} \cdot 163) \\
a_{7} &= 2^{2} \cdot (j - 2^{6} \cdot 3^{3})^{3} \cdot j^{4} \cdot (11 \cdot 61 \cdot 193 \cdot j^{5} - 2^{3} \cdot 3 \cdot 27510443 \cdot j^{4} \\
&\quad + 2^{9} \cdot 3^{3} \cdot 97550587 \cdot j^{3} - 2^{16} \cdot 3^{6} \cdot 11 \cdot 2599451 \cdot j^{2} \\
\displaybreak[0]
&\quad + 2^{23} \cdot 3^{9} \cdot 5 \cdot 739 \cdot 1109 \cdot j - 2^{34} \cdot 3^{13} \cdot 4691) \\
a_{6} &= 2^{3} \cdot (j - 2^{6} \cdot 3^{3})^{3} \cdot j^{4} \cdot (2^{4} \cdot 3^{2} \cdot j^{8} + 7 \cdot 199 \cdot 1373 \cdot j^{7} \\
&\qquad - 2^{2} \cdot 29 \cdot 37 \cdot 281 \cdot 13913 \cdot j^{6} + 2^{13} \cdot 3^{3} \cdot 7 \cdot 233 \cdot 143281 \cdot j^{5} \\
&\quad - 2^{15} \cdot 3^{7} \cdot 5 \cdot 11 \cdot 21117827 \cdot j^{4} + 2^{23} \cdot 3^{9} \cdot 3943 \cdot 117577 \cdot j^{3} \\
\displaybreak[0]
&\quad - 2^{31} \cdot 3^{12} \cdot 769 \cdot 45317 \cdot j^{2} + 2^{41} \cdot 3^{16} \cdot 7 \cdot 15923 \cdot j - 2^{50} \cdot 3^{20} \cdot 269) \\
a_{5} &= 2^{4} \cdot (j - 2^{6} \cdot 3^{3})^{4} \cdot j^{5} \cdot (2^{6} \cdot 3^{4} \cdot 5 \cdot j^{8} - 7 \cdot 5051 \cdot 5939 \cdot j^{7} \\
&\quad + 2^{3} \cdot 3^{2} \cdot 5 \cdot 61 \cdot 101 \cdot 330037 \cdot j^{6} - 2^{9} \cdot 3^{5} \cdot 96289 \cdot 119173 \cdot j^{5} \\
&\quad + 2^{16} \cdot 3^{9} \cdot 17 \cdot 77252741 \cdot j^{4} - 2^{22} \cdot 3^{11} \cdot 11 \cdot 71 \cdot 523 \cdot 4091 \cdot j^{3} \\
\displaybreak[0]
&\quad + 2^{35} \cdot 3^{14} \cdot 5 \cdot 673 \cdot 977 \cdot j^{2} - 2^{41} \cdot 3^{18} \cdot 79 \cdot 1831 \cdot j + 2^{55} \cdot 3^{24}) \\
a_{4} &= (j - 2^{6} \cdot 3^{3})^{4} \cdot j^{6} \cdot (2^{8} \cdot 3^{3} \cdot 5 \cdot 2003 \cdot j^{9} - 409 \cdot 39157 \cdot 44483 \cdot j^{8} \\
&\quad + 2^{9} \cdot 3 \cdot 2092618568983 \cdot j^{7} - 2^{20} \cdot 3^{4} \cdot 98512996093 \cdot j^{6} \\
&\quad + 2^{20} \cdot 3^{7} \cdot 41 \cdot 242261 \cdot 608831 \cdot j^{5} - 2^{28} \cdot 3^{10} \cdot 5 \cdot 1231 \cdot 155631757 \cdot j^{4} \\
&\quad + 2^{32} \cdot 3^{13} \cdot 521 \cdot 3077579657 \cdot j^{3} - 2^{42} \cdot 3^{16} \cdot 997 \cdot 1607 \cdot 16657 \cdot j^{2} \\
\displaybreak[0]
&\quad + 2^{52} \cdot 3^{20} \cdot 23 \cdot 541 \cdot 6863 \cdot j - 2^{63} \cdot 3^{24} \cdot 5 \cdot 11987) \\
a_{3} &= 2 \cdot (j - 2^{6} \cdot 3^{3})^{5} \cdot j^{6} \cdot (3^{2} \cdot 377732207 \cdot j^{10} - 2^{6} \cdot 5^{2} \cdot 7 \cdot 101 \cdot 28520381 \cdot j^{9} \\
&\quad + 2^{11} \cdot 11 \cdot 337 \cdot 17990477821 \cdot j^{8} - 2^{20} \cdot 3^{3} \cdot 179 \cdot 389 \cdot 171956657 \cdot j^{7} \\
&\quad + 2^{23} \cdot 3^{6} \cdot 5 \cdot 479 \cdot 37193046587 \cdot j^{6} - 2^{30} \cdot 3^{9} \cdot 1283 \cdot 28703 \cdot 758137 \cdot j^{5} \\
&\quad + 2^{36} \cdot 3^{12} \cdot 7 \cdot 31 \cdot 54791988203 \cdot j^{4} - 2^{45} \cdot 3^{15} \cdot 19^{2} \cdot 151 \cdot 7738067 \cdot j^{3} \\
\displaybreak[0]
&\quad + 2^{55} \cdot 3^{20} \cdot 41 \cdot 12810583 \cdot j^{2} - 2^{65} \cdot 3^{24} \cdot 1103107 \cdot j + 2^{76} \cdot 3^{27} \cdot 1447) \\
a_{2} &= 2^{2} \cdot (j - 2^{6} \cdot 3^{3})^{5} \cdot j^{7} \cdot (42967 \cdot 2406947 \cdot j^{11} - 2^{3} \cdot 557 \cdot 1783 \cdot 140768209 \cdot j^{10} \\
&\quad + 2^{9} \cdot 3^{4} \cdot 6205891 \cdot 21226039 \cdot j^{9} - 2^{19} \cdot 3^{7} \cdot 5 \cdot 11 \cdot 251872948013 \cdot j^{8} \\
&\quad + 2^{24} \cdot 3^{9} \cdot 5 \cdot 13 \cdot 23 \cdot 37 \cdot 521 \cdot 3203149 \cdot j^{7} - 2^{29} \cdot 3^{13} \cdot 47242981376477 \cdot j^{6} \\
&\quad + 2^{35} \cdot 3^{16} \cdot 227 \cdot 112292655271 \cdot j^{5} - 2^{41} \cdot 3^{18} \cdot 107 \cdot 269749728667 \cdot j^{4} \\
&\quad + 2^{54} \cdot 3^{22} \cdot 43 \cdot 449215127 \cdot j^{3} - 2^{61} \cdot 3^{27} \cdot 5 \cdot 653 \cdot 54193 \cdot j^{2} \\
\displaybreak[0]
&\quad + 2^{72} \cdot 3^{30} \cdot 139 \cdot 3719 \cdot j - 2^{82} \cdot 3^{35} \cdot 139) \\
a_{1} &= 2^{3} \cdot (j - 2^{6} \cdot 3^{3})^{6} \cdot j^{8} \cdot (1847032397279 \cdot j^{11} - 2^{6} \cdot 47 \cdot 157 \cdot 3691 \cdot 11660843 \cdot j^{10} \\
&\quad + 2^{14} \cdot 3^{4} \cdot 383 \cdot 25679 \cdot 7797631 \cdot j^{9} - 2^{20} \cdot 3^{6} \cdot 400129001343469 \cdot j^{8} \\
&\quad + 2^{24} \cdot 3^{9} \cdot 5 \cdot 41 \cdot 503 \cdot 67307 \cdot 267271 \cdot j^{7} \\
&\quad - 2^{30} \cdot 3^{12} \cdot 19 \cdot 509 \cdot 13597 \cdot 11431571 \cdot j^{6} \\
&\quad + 2^{37} \cdot 3^{15} \cdot 31 \cdot 3038701 \cdot 4610147 \cdot j^{5} - 2^{43} \cdot 3^{20} \cdot 7^{2} \cdot 41 \cdot 73 \cdot 2381 \cdot 56891 \cdot j^{4} \\
&\quad + 2^{52} \cdot 3^{21} \cdot 5 \cdot 139 \cdot 9239401667 \cdot j^{3} - 2^{62} \cdot 3^{25} \cdot 5 \cdot 1381 \cdot 3698087 \cdot j^{2} \\
\displaybreak[0]
&\quad + 2^{73} \cdot 3^{29} \cdot 11 \cdot 47 \cdot 58693 \cdot j - 2^{85} \cdot 3^{33} \cdot 8161) \\
a_{0} &= - 2^{4} \cdot (j - 2^{6} \cdot 3^{3})^{6} \cdot j^{8} \cdot (2^{3} \cdot 3^{2} \cdot 7^{6} \cdot j^{14} - 5 \cdot 13 \cdot 3109 \cdot 76441597 \cdot j^{13} \\
&\quad + 2^{4} \cdot 3449 \cdot 4363 \cdot 873750089 \cdot j^{12} - 2^{11} \cdot 3^{4} \cdot 7 \cdot 2087 \cdot 57859 \cdot 9420337 \cdot j^{11} \\
&\quad + 2^{16} \cdot 3^{8} \cdot 11^{2} \cdot 73 \cdot 125183 \cdot 10636957 \cdot j^{10} - 2^{26} \cdot 3^{9} \cdot 691 \cdot 14434308694753 \cdot j^{9} \\
&\quad + 2^{31} \cdot 3^{13} \cdot 101 \cdot 283 \cdot 252059913139 \cdot j^{8} \\
&\quad - 2^{37} \cdot 3^{16} \cdot 11 \cdot 13 \cdot 17 \cdot 647 \cdot 863 \cdot 4253233 \cdot j^{7} \\
&\quad + 2^{43} \cdot 3^{18} \cdot 631819 \cdot 16451871913 \cdot j^{6} - 2^{48} \cdot 3^{23} \cdot 149 \cdot 233 \cdot 90533 \cdot 330413 \cdot j^{5} \\
&\quad + 2^{59} \cdot 3^{25} \cdot 23 \cdot 1408302006413 \cdot j^{4} - 2^{70} \cdot 3^{27} \cdot 726838208711 \cdot j^{3} \\
\displaybreak[0]
&\quad + 2^{80} \cdot 3^{32} \cdot 7 \cdot 263 \cdot 337 \cdot 1327 \cdot j^{2} - 2^{90} \cdot 3^{37} \cdot 569731 \cdot j + 2^{100} \cdot 3^{39} \cdot 17^{3}) \\
\displaybreak[0]
b_{11} &= - (j - 2^{6} \cdot 3^{3}) \cdot j \\
\displaybreak[0]
b_{10} &= - 2 \cdot 13 \cdot 3^{2} \cdot (j - 2^{6} \cdot 3^{3}) \cdot j^{2} \\
\displaybreak[0]
b_{9} &= 2^{2} \cdot (j - 2^{3} \cdot 3^{6}) \cdot (j - 2^{6} \cdot 3^{3})^{2} \cdot j^{2} \\
\displaybreak[0]
b_{8} &= 3^{4} \cdot (13 \cdot j - 2^{5} \cdot 3 \cdot 163) \cdot (j - 2^{6} \cdot 3^{3})^{2} \cdot j^{3} \\
\displaybreak[0]
b_{7} &= 5 \cdot 2^{5} \cdot 3^{6} \cdot (j - 2^{6} \cdot 3^{3})^{3} \cdot j^{4} \\
\displaybreak[0]
b_{6} &= 2^{2} \cdot 3^{3} \cdot (j - 2^{6} \cdot 3^{3})^{3} \cdot j^{4} \cdot (j^{2} + 2^{4} \cdot 3^{5} \cdot 13 \cdot j - 2^{9} \cdot 3^{5} \cdot 269) \\
\displaybreak[0]
b_{5} &= 2^{5} \cdot 3^{5} \cdot (5 \cdot j - 2^{6} \cdot 3^{4}) \cdot (j - 2^{6} \cdot 3^{3})^{4} \cdot j^{5} \\
\displaybreak[0]
b_{4} &= 2^{5} \cdot 3^{8} \cdot (31 \cdot j - 2^{3} \cdot 3^{2} \cdot 1471) \cdot (j - 2^{6} \cdot 3^{3})^{4} \cdot j^{6} \\
\displaybreak[0]
b_{3} &= 2^{8} \cdot 3^{8} \cdot (383 \cdot j - 2^{6} \cdot 3 \cdot 1447) \cdot (j - 2^{6} \cdot 3^{3})^{5} \cdot j^{6} \\
\displaybreak[0]
b_{2} &= 2^{9} \cdot 3^{9} \cdot (3923 \cdot j - 2^{6} \cdot 3^{5} \cdot 139) \cdot (j - 2^{6} \cdot 3^{3})^{5} \cdot j^{7} \\
\displaybreak[0]
b_{1} &= 13 \cdot 19 \cdot 3^{11} \cdot 2^{15} \cdot (j - 2^{6} \cdot 3^{3})^{6} \cdot j^{8} \\
b_{0} &= - 2^{8} \cdot 3^{9} \cdot (j - 2^{6} \cdot 3^{3})^{6} \cdot j^{8} \cdot (j^{2} - 2^{7} \cdot 3^{3} \cdot 1399 \cdot j + 2^{12} \cdot 3^{6} \cdot 17^{3})
\end{align*}

\end{document}